\documentclass[letterpaper, 10 pt, conference]{ieeeconf}  % Comment this line out

\usepackage{amsmath}
\usepackage{amssymb}
\usepackage{amsfonts}
\usepackage{graphicx}
\usepackage{subfigure}
\usepackage{url}
\usepackage{color}

\newcommand{\change}[1]{#1}

\newcommand{\note}[1]{\textcolor{blue}{{#1}}}

\newcommand{\trace}{{\rm Tr}}
\newcommand{\diag}{{\rm Diag}}
\newcommand{\rank}{{\rm rank}}

\newcommand{\FixedRankPSD}[2]{S_+({#1},{#2})}

\newcommand{\EDM}[1]{\mathrm{EDM}({#1})}

\newcommand{\OG}[1]{{\mathcal{O}({#1})}}

\newcommand{\set}[2]{{\{{#1}\,:\ {#2}\}}}

\newcommand{\grad}{\mathrm{grad}}
\newcommand{\hess}{\mathrm{Hess}}

\newcommand{\mat}[1]{{\bf #1}}
\renewcommand{\vec}[1]{{\bf #1}}

\IEEEoverridecommandlockouts                              % This command is only
\title{\LARGE \bf
Low-rank optimization for distance matrix completion}

\author{B. Mishra, G. Meyer and R. Sepulchre% <-this % stops a space
\thanks{The authors are with the Department of Electrical Engineering and Computer Science, University of Li\`ege, Montefiore Institute, Sart-Tilman, 4000 Li\`ege, Belgium (e-mails: \{b.mishra, g.meyer, r.sepulchre\}@ulg.ac.be).}
}

\begin{document}

\maketitle
\thispagestyle{empty}
\pagestyle{empty}

\begin{abstract}
This paper addresses the problem of low-rank distance matrix completion. This problem amounts to recover the missing entries of a distance matrix when the dimension of the data embedding space is possibly unknown but small compared to the number of considered data points. The focus is on high-dimensional problems. We recast the considered problem into an optimization problem over the set of low-rank positive semidefinite matrices and propose two efficient algorithms for low-rank distance matrix completion. In addition, we propose a strategy to determine the dimension of the embedding space. The resulting algorithms scale to high-dimensional problems and monotonically converge to a global solution of the problem. Finally, numerical experiments illustrate the good performance of the proposed algorithms on benchmarks.

\vspace{1em}
 \note{\emph{This is the pre-print version of \cite{mishra11a}.}}

\end{abstract}

\section{INTRODUCTION}
Completing the missing entries of a matrix under low-rank constraint is a fundamental and recurrent problem in many modern engineering applications (see \cite{candes09a} and references therein). Recently, the problem has gained much popularity thanks to collaborative filtering applications and the Netflix challenge \cite{netflix09a}.

This paper focuses on an important variant of the problem, that is, completing the missing entries of a Euclidean distance matrix (EDM) under low-rank constraint. Typical applications include data visualization \cite{morrison03a}, dimensionality reduction in behavioral sciences and economics \cite{borg97a}, molecular conformation problems \cite{more99a, glunt93a}, just to name a few.

A Euclidean distance matrix $\mat{D}\in\mathbb{R}^{n\times n}$ contains the (squared) pairwise distances between $n$ data points $\vec{y}_i\in\mathbb{R}^{r}$, $i=1,...,n$. This matrix is symmetric and has a zero diagonal. Its entries are non-negative and satisfy the triangle inequality. These properties are readily verified by examining the entries of the distance matrix,
\begin{equation*}
	\mat{D}_{ij}=\|\vec{y}_i - \vec{y}_j\|_2^2.
\end{equation*}
The set $\EDM{n}$ of $n$-by-$n$ Euclidean distance matrices forms a convex cone which has a well-studied geometry (see \cite{dattorro05a, alfakih05a}, and references therein). One property of a Euclidean distance matrix is that it is rank deficient. The rank of $\mat{D}$ is upper bounded by $r+2$ (and the rank is generically $r+2$), which in many problems is very small compared to $n$, the number of data points.

Given a set of pairwise distances or dissimilarities between data points, the goal of low-rank distance matrix completion algorithms is to recover a full Euclidean distance matrix from a restrictive set of given distances. Inference on the unknown entries is possible thanks to the low-rank property which models the redundancy between the available data.

A closely related problem is multidimensional scaling (MDS) for which all pairwise distances are available up front. A solution to this problem is the classical multidimensional scaling algorithm (CMDS), which relies on singular value decomposition to find a globally optimum embedding of fixed-rank. The CMDS algorithm minimizes the total quadratic error on scalar products between data points. Other algorithms have focused on variant cost functions, see the paper \cite{de-leeuw01a} for a survey in this area.

In contrast to the classical multidimensional scaling formulation, the problem of Euclidean distance matrix completion involves missing distances. The problem can be considered as a variant of multidimensional scaling problem with binary weights \cite{de-leeuw01a,kearsley98a}. The low-rank distance matrix completion problem is known to be NP-hard in general \cite{laurent98a, huang03a}, but convex relaxations have been proposed to render the problem  tractable \cite{alfakih99a, cayton06a}. Typical convex relaxations cast the EDM completion problem into a convex optimization problem on the set of positive semidefinite matrix, resulting in semidefinite programming techniques \cite{boyd04a}. This convex formulation is nevertheless a large-scale problem when $n$ is large.

Imposing the rank constraint in the problem formulation is an appealing way of reducing the size of the search space. However, it results in a non-convex optimization problem. Although convergence results are only local, the approach performs well in practice \cite{trosset00a}. Both first-oder \cite{kruskal64a,buja04a} and second order \cite{glunt93a,tarazaga93a,kearsley98a,fang10a} optimization methods have been considered and heuristics for finding a good low-rank initialization have been proposed \cite{fang10a}. 

A difficulty encountered by second order optimization algorithms is the intrinsic invariance properties of the data representation due to rotations. This issue may indeed prevent second order optimization algorithms to converge \cite{absil09a}. Several authors have resolved this issue at the extra cost of normalizing the data representation \cite{kearsley98a} or adding a penalization term to the objective function \cite{tarazaga93a}. In this paper, the invariance to rotations is lifted in the problem formulation and is free of additional computational cost (see Section \ref{sec:manifold}). A survey of low-rank distance matrix completion algorithms can be found in the recent papers \cite{chu03a,fang10a}. 

Although, the problem is not new and is well-studied, a practical limitation of most of existing algorithms is that they do not scale to high-dimensional problems. Moreover, the problem of choosing a priori an appropriate dimension for the data embedding is still an open research question. 

In this paper, the focus is on efficient algorithms that scale to high-dimensional problems. Following a number of previous contributions in the literature, we recast the problem into an optimization problem over the set of low-rank positive semidefinite matrices. We adopt the geometric optimization framework of optimization on Riemannian matrix manifolds \cite{absil08a}. \change{Our main contribution is to extend the framework developed in \cite{journee10a} to the problem of low-rank distance matrix completion. This results in an efficient strategy for estimating the dimension of the embedding space. The proposed algorithms have linear complexity in the problem size and in the number of available distances. The strategy for estimating the optimal embedding dimension ensures that the proposed algorithms converge monotonically to the global (low-rank) solution of the problem.}

The paper is organized as follow. Section \ref{sec:formulation} presents the problem of interest and its different formulations. Section \ref{sec:manifold} describes the chosen optimization framework and introduces the main geometrical objects required by our algorithms. Section \ref{sec:algorithms} is devoted to the design of efficient algorithms for low-rank distance matrix completion. Finally, Section \ref{sec:experiments} presents some numerical simulations.

\section{LOW-RANK DISTANCE MATRIX COMPLETION}\label{sec:formulation}
Given a set of dissimilarities $\widetilde{\mat{D}}_{ij} \geq 0$ between $n$ data points, distance matrix completion algorithms solve
\begin{equation}\label{eq:formulation-on-edm}
	\min_{\mat{D}\in\EDM{n}} \|\mat{H}\odot(\mat{D}-\widetilde{\mat{D}})\|_F^2,
\end{equation}
where $\mat{H}$ is a symmetric matrix with binary entries and the operator $\odot$ denotes elementwise multiplication. If $\mathcal{D}$ is the set of given entries $(i,j)$ in $\widetilde{\mat{D}}$ such that $i < j$, then
\begin{equation*}
	\mat{H}_{ij} = \mat{H}_{ji} =
	\begin{cases} 1 & \text{if $(i,j)\in\mathcal{D}$,}\\
	0 &\text{otherwise.}
	\end{cases}
\end{equation*}
The number of elements in the set $\mathcal{D}$ is denoted by $d$. Although, $d$ is at most equal to $n(n-1)/2$, in most applications, it is of order $O(nr)$, where $r$ is the optimal embedding dimension. Dissimilarities potentially differ from distances in that they are not required to satisfy triangle inequality. For instance, this takes into account the fact that observation noise could make $\widetilde{\mat{D}}$ different from a valid EDM.

A convenient alternative formulation of \eqref{eq:formulation-on-edm} is to cast this problem into an optimization problem on the set of positive semidefinite matrices \cite{alfakih99a}. \change{The reformulation hinges on a classical result by Schoenberg which relates Euclidean distance matrices and positive semidefinite matrices of rank equal to the dimension of the embedding space \cite{schoenberg36a}. The corresponding reformulation can be written as}
\begin{equation}\label{eq:formulation-on-psd}
	\min_{\mat{X}\succeq 0}\; \|\mat{H}\odot(\kappa(\mat{X})-\widetilde{\mat{D}})\|_F^2,
\end{equation}
where $\kappa$ is a mapping from the set of positive semidefinite matrices to the set of Euclidean distance matrices
\begin{equation*}
	\kappa(\mat{X})=\diag(\mat{X})\vec{1}^T + \vec{1}\diag(\mat{X})^T - 2\mat{X}.
\end{equation*}
The function $\diag(\cdot)$ extracts the diagonal of its argument, and $\vec{1}$ denotes a vector with all entries equal to one.

A practical advantage of \eqref{eq:formulation-on-psd} compared to \eqref{eq:formulation-on-edm} is that the rank of $\mat{X}$ identifies with the dimension of the embedding space. When no restriction is imposed on the rank of $\mat{X}$, problem \eqref{eq:formulation-on-psd} is convex and thus presents a global solution. 

In this paper, we consider the case where the global solution $\mat{X}^{\star}$ of \eqref{eq:formulation-on-psd} is low-rank that is,
\begin{equation}\label{eq:low-rank-assumption}
	\rank(\mat{X}^{\star}) = r \ll n.
\end{equation}
Following \cite{journee10a}, we solve a sequence of nonconvex problems of increasing dimension until the actual value of the rank $r$ is reached. Each nonconvex problem consists in solving the following rank-constrained optimization problem
\begin{equation}\label{eq:formulation-on-low-rank-psd}
	\min_{\mat{X}\succeq 0}\ \|\mat{H}\odot(\kappa(\mat{X})-\widetilde{\mat{D}})\|_F^2,\quad \text{s.t.}\quad \rank(\mat{X}) = p.
\end{equation}

By screening values from $p=1$ to $p=r$, the results presented in \cite{journee10a} guarantee a monotonic convergence to a solution of the original problem \eqref{eq:formulation-on-psd}. The proposed strategy for finding the actual rank $r$ is detailed in Section \ref{sec:strategy-rank}.

Problem \eqref{eq:formulation-on-low-rank-psd} is solved efficiently by exploiting a low-rank parametrization of the search space. The proposed approach hinges on the fact that any rank-$p$ positive semidefinite matrix admits a factorization
\begin{equation*}\label{eq:low-rank-factorization}
	\mat{X} = \mat{Y}\mat{Y}^T,
\end{equation*}
where $\mat{Y}\in\mathbb{R}_*^{n\times p}= \set{\mat{Y}\in\mathbb{R}^{n\times p}}{\det(\mat{Y}^T\mat{Y})\neq 0}$. 

To exploit this factorization, we adopt the geometric framework of optimization on Riemannian manifolds \cite{absil08a}. Basic concepts and notations are introduced in the next section. See the book \cite{absil08a} for further details on optimization on matrix manifolds and for a state-of-the-art in this area.

\section{MANIFOLD-BASED OPTIMIZATION}\label{sec:manifold}
An intrinsic property of the factorization $\mat{X}=\mat{Y}\mat{Y}^T$ is that it is invariant with respect to the transformation
\begin{equation*}
	\mat{Y}\mapsto \mat{Y}\mat{Q},
\end{equation*}
where $\mat{Q}\in\OG{p}=\set{\mat{Q}\in\mathbb{R}^{p\times p}}{\mat{Q}^T\mat{Q}=\mat{Q}\mat{Q}^T=\mat{I}}$.

This invariance property renders the minima of a cost function $f(\mat{Y}\mat{Y}^T)$ not isolated. This issue is not harmful for first order-methods such as gradient descent algorithms but greatly affects the convergence properties of second-order methods \cite{absil08a, absil09a}. 

To circumvent this issue, we reformulate the problem of interest as an optimization problem on the quotient manifold
\begin{equation}\label{eq:search-space}
	\mathcal{M}\triangleq \FixedRankPSD{p}{n}\simeq \mathbb{R}_*^{n\times p}/\OG{p},
\end{equation}
which represents the set of equivalence classes
\begin{equation}\label{eq:equiv-classes}
	[\mat{Y}] =\set{\mat{Y}\mat{Q}}{\mat{Q}\in \OG{p}}.
\end{equation}
The set $\FixedRankPSD{p}{n}$ is the set of rank-$p$ symmetric positive semidefinite matrices of size $n$, that is,
\begin{equation*}
	\FixedRankPSD{p}{n} = \set{\mat{X}\in\mathbb{R}^{n\times n}}{\mat{X}=\mat{X}^T\succeq 0,\;\rank(\mat{X})=p}.
\end{equation*}
This set has a rich Riemannian manifold geometry which can be exploited for algorithmic purposes \cite{vandereycken09a, bonnabel09a, meyer10a}.

Problem \eqref{eq:formulation-on-low-rank-psd} is now reformulated as an unconstrained optimization problem over the set of equivalence classes \eqref{eq:equiv-classes},
\begin{equation}\label{eq:formulation-manifold}
	\min_{[\mat{Y}]\in\mathcal{M}} f([\mat{Y}]),
\end{equation} 
for the cost function 
\begin{equation}\label{eq:cost-function}
	f([\mat{Y}]) = \|\mat{H}\odot(\kappa(\mat{Y}\mat{Y}^T)-\widetilde{\mat{D}})\|_F^2.
\end{equation}
To develop optimization algorithms on the quotient manifold, \change{the tangent space $T_{\mat{Y}}\mathcal{M}$ of \eqref{eq:search-space}} is endowed with the Riemannian metric
\begin{equation*}
	g_{\mat{Y}}(\xi_{\mat{Y}},\eta_{\mat{Y}}) = \trace(\xi_{\mat{Y}}^{\;T}\eta_{\mat{Y}}),\quad \xi_{\mat{Y}},\eta_{\mat{Y}}\in T_{\mat{Y}}\mathcal{M},
\end{equation*}
which is inherited from the natural metric of $\mathbb{R}^{n\times p}$. With this metric, the tangent space $T_{\mat{Y}}\mathcal{M}$ at a given point $\mat{Y}$ is decomposed into the sum of two complementary spaces,
\begin{equation*}
	T_{\mat{Y}}\mathcal{M} = \mathcal{V}_{\mat{Y}}\mathcal{M}\oplus\mathcal{H}_{\mat{Y}}\mathcal{M}.
\end{equation*}
The vertical space $\mathcal{V}_{\mat{Y}}\mathcal{M}$ contains the set of directions that are tangent to the set of equivalence classes \eqref{eq:equiv-classes}, that is, 
\begin{equation*}
	\mathcal{V}_{\mat{Y}}\mathcal{M} = \set{\mat{Y}\mat{\Omega}}{\mat{\Omega}^T=-\mat{\Omega}\in\mathbb{R}^{p\times p}}.
\end{equation*}
The horizontal space $\mathcal{H}_{\mat{Y}}\mathcal{M}$ contains the directions $\bar{\xi}_{\mat{Y}}$ that are orthogonal to the set of equivalence classes, 
\begin{equation*}
	\mathcal{H}_{\mat{Y}}\mathcal{M} = \set{\bar{\xi}_{\mat{Y}}\in\mathbb{R}^{n\times p}}{\bar{\xi}_{\mat{Y}}^{\,T}\mat{Y} = \mat{Y}^T\bar{\xi}_{\mat{Y}}}.
\end{equation*}
With such a construction, the directions of interest can be restricted to horizontal directions $\bar{\xi}_{\mat{Y}}$. Indeed, displacements along vertical directions leave the cost function unchanged.

The projection of a direction $\mat{Z}\in\mathbb{R}^{n\times p}$ onto the horizontal space is given by $\Pi_{\mathcal{H}_{\mat{Y}}}(\mat{Z}) = \mat{Z} - \mat{Y}\mat{\Omega}$, where $\mat{\Omega}\in\mathbb{R}^{p\times p}$ is skew-symmetric and satisfies the Sylvester equation
\begin{equation*}
	\mat{\Omega Y}^T \mat{Y} + \mat{Y}^T \mat{Y \Omega} =\mat{Y}^T\mat{ Z} -\mat{ Z}^T \mat{Y}.
\end{equation*}
Overall, projecting a direction $\mat{Z}\in\mathbb{R}^{n\times p}$ onto the horizontal space requires $O(np^2 + np + p^3)$ operations (computing matrices $\mat{Y}^T \mat{Y}$, $\mat{Y}^T\mat{Z}$, and $ \mat{Y}\mat{\Omega}$ requires $O(np^2)$ operations, solving the Sylvester equation is performed in $O(p^3)$ operations and the projection requires $O(np)$ operations).

To update our search variable, we require a local mapping from tangent space to the manifold. Such a mapping is called a retraction. For the manifold of interest, a retraction is provided by the simple and efficient formula
\begin{equation}\label{eq:retraction}
	\mathrm{R}_{\mat{Y}}(\bar{\xi}_\mat{Y}) = \mat{Y} + \bar{\xi}_\mat{Y}.
\end{equation}
which gives a full-rank matrix for generic direction $\bar{\xi}_\mat{Y}$.
\section{ALGORITHMS}\label{sec:algorithms}
In this section, we exploit the concepts presented in the previous section to develop both a gradient descent algorithm and a trust-region algorithm to solve \eqref{eq:formulation-manifold}.

\subsection{Gradient descent algorithm}
The gradient of a smooth cost function $f:\mathcal{M}\rightarrow\mathbb{R}$ is the unique tangent vector $\grad f(\mat{Y})\in T_\mat{Y}\mathcal{M}$ that satisfies
\begin{equation}\label{eq:gradient-formula}
	g_{\mat{Y}}(\xi_{\mat{Y}},\grad f(\mat{Y})) = Df(\mat{Y})[\xi_{\mat{Y}}],\quad \forall \xi_{\mat{Y}}\in T_\mat{Y}\mathcal{M}.
\end{equation}
The quantity $Df(\mat{Y})[\xi_{\mat{Y}}]$ is the directional derivative of $f$ in the direction $\xi_{\mat{Y}}$, that is,
\begin{equation*}
	Df(\mat{Y})[\xi_{\mat{Y}}] = \lim_{t\rightarrow 0}\frac{f(\mat{Y} + t\xi_{\mat{Y}})-f(\mat{Y})}{t}.
\end{equation*}
Applying formula \eqref{eq:gradient-formula} to the cost \eqref{eq:cost-function} gives us the gradient
\begin{equation}\label{eq:gradient-cost}
	\grad f(\mat{Y}) = 2\kappa^*(\mat{H}\odot(\kappa(\mat{YY}^T )-\widetilde{\mat{D}}))\mat{Y},
\end{equation}
where $\kappa^*(\mat{A})$ is the adjoint operator of $\kappa$ defined by
\begin{equation*}
\kappa^*(\mat{A}) = 2(\diag(\mat{A}\vec{1}) - \mat{A}).
\end{equation*}
Combining the gradient \eqref{eq:gradient-cost} with the retraction \eqref{eq:retraction} gives us the gradient descent algorithm
\begin{equation}\label{eq:gradient-descent-algo}
	\mat{Y}_{t+1} = \mat{Y}_t - 2s_t\kappa^*(\mat{H}\odot(\kappa(\mat{Y}_t\mat{Y}_t^T )-\widetilde{\mat{D}}))\mat{Y}_t,
\end{equation}
where $s_t > 0$ is the gradient step size. We select $s_t$ using the Armijo criterion \cite{nocedal06a}, that is, a step size $s_A$ that satisfies
\begin{equation*}
	f(\mat{Y}_{t} - s_A\,\grad f(\mat{Y}_t)) \leq f(\mat{Y}_{t}) - c\,s_A \|\grad f(\mat{Y}_t)\|_F^2,
\end{equation*}
where $c\in(0,1)$ is a constant (we choose the value $c=0.5$).

The asymptotic computational cost of an iteration \eqref{eq:gradient-descent-algo} is $O( d p + np)$, where $d$ is the number of known entries of $\widetilde{\mat{D}}$. The memory requirement is \change{$O(d+np)$}. The computationally most demanding step is the computation of the gradient, which requires $O(dp)$ operations. This low computational complexity and memory requirement allows us to handle potentially large data sets. A drawback is however that the gradient descent algorithm only guarantees a linear convergence rate. We can achieve a superlinear convergence rate by means of a Riemannian trust-region algorithm which exploits second-order information.
%Since $d$ is typically of order $O(np)$, computing the gradient thus typically requires $O(np^2)$ operations. 

\subsection{Trust-region algorithm}
Trust-region algorithms sequentially solve the problem
\begin{align*}
	\min_{\bar{\xi}\in\mathcal{H}_{\mat{Y}}\mathcal{M}} 
	&f(\mat{Y}) + g_{\mat{Y}}(\bar{\xi},\grad\,f(\mat{Y})) 
	+ \frac{1}{2}g_{\mat{Y}}(\bar{\xi},\hess\,f(\mat{Y})[\bar{\xi}]),\\
	& \text{s.t.}\quad g_{\mat{Y}}(\bar{\xi},\bar{\xi})\leq\delta^2,
\end{align*}
which amounts to minimize a quadratic model of the cost function on a trust-region radius of size $\delta$. Once a search direction $\bar{\xi}$ is identified, the search variable is updated as
\begin{equation}
	\mat{Y}_{t+1} = \mathrm{R}_{\mat{Y}_t}(\bar{\xi}).
\end{equation}

The trust-region radius $\delta$ vary according to the quality of the iterate. When a good solution is found within the trust-region, then the trust-region is expanded. Conversely, if the iterate is poor then the region is contracted.

More technical details on trust-region algorithms on Riemannian manifolds can be found in \cite{absil07a, absil08a}. In this paper, we adapt the generic implementation of the toolbox GenRTR to our problem of interest.\footnote{The software can be downloaded from\\\url{http://www.math.fsu.edu/~cbaker/GenRTR/}}

Trust-region algorithms require the computation of the Riemannian Hessian $\hess\,f(\mat{Y})[\bar{\eta}]$ in a given direction $\bar{\eta}$. It is obtained as
\begin{equation*}
	\hess\,f(\mat{Y})[\bar{\eta}] 
		\triangleq \nabla_{\bar{\eta}}\;\grad\,f(\mat{Y})
		%= \Pi_{\mathcal{H}_{\mat{Y}}}(D\grad\,f(\mat{Y})[\bar{\eta}])
\end{equation*}
where $\nabla_{\bar{\eta}}\,\grad\,f(\mat{Y})$ is the Riemannian connection of the gradient vector field in the direction $\bar{\eta}$. Riemannian connections generalize the notion of directional derivative of a vector field to Riemannian manifolds. Given a vector field $\zeta$ on $\mathcal{M}$ that assigns to each point $\mat{Y}\in\mathcal{M}$ a tangent vector $\zeta_{\mat{Y}}\in T_{\mat{Y}}\mathcal{M}$, the directional derivative of $\zeta$ at $\mat{Y}\in\mathcal{M}$ in a direction $\bar{\eta}\in\mathcal{H}_{\mat{Y}}\mathcal{M}$ is given by
\begin{equation}
	\nabla_{\bar{\eta}}\;\zeta_{\mat{Y}} 
	= \Pi_{\mathcal{H}_{\mat{Y}}}\left(\lim_{t\rightarrow 0}\frac{\zeta_{\mat{Y} + t\bar{\eta}} - \zeta_{\mat{Y}}}{t}\right).
\end{equation}
Applying this formula to the vector field $\grad\,f(\mat{Y})$ gives us
\begin{align*}
	\hess\,f(\mat{Y})[\bar{\eta}] = 2\Pi_{\mathcal{H}_{\mat{Y}}}(&\kappa^*(\mat{H}\odot(\kappa(\mat{Y}\bar{\eta}^{\,T} + \bar{\eta} \mat{Y}^T)))\mat{Y}\\
	&+\kappa^*(\mat{H}\odot (\kappa(\mat{Y}\mat{Y}^T)- \widetilde{\mat{D}}))\bar{\eta}).
\end{align*}

The numerical cost of an iteration of the trust-region algorithm is $O(dp + np + np^2 + p^3)$. The memory requirement is \change{$O(d + np )$}. The computational bottleneck is the computation of the Hessian. Still, the complexity is linear in both the number of available distance and in the problem size. With a proper parameter tuning, the proposed trust-region algorithm enjoys a superlinear convergence rate.

\subsection{Strategy for estimating the optimal embedding dimension}\label{sec:strategy-rank}
The following section is an adaptation of the material presented in \cite{journee10a} to the problem of interest. To identify the (unknown) rank $r$ of the global solution to \eqref{eq:formulation-on-psd}, we solve a sequence of nonconvex problems \eqref{eq:formulation-manifold} of increasing dimension. The approximation rank $p$ is progressively incremented from $p=1$ to $p=r$. Using a warm restart strategy for moving from one value of $p$ to the next, we are able to propose a descent algorithm that converges monotonically to a global solution of  the original problem \eqref{eq:formulation-on-low-rank-psd}. 

This strategy efficiently exploits the previous iterations of the algorithm as opposed to earlier heuristic methods that use random restart for each value of the rank \cite{burer03a}.

For a given rank $p < r$, the trust-region or gradient descent algorithm gives us a local minimizer $\mat{Y}_{p}^{\star}$ of the nonconvex problem \eqref{eq:formulation-manifold}. Let us consider the following initial condition for the problem of rank $p+1$,
\begin{equation*}
	\mat{Y}_{0} = [\mat{Y}_{p}^{\star}|\mat{0}^{n\times 1}],
\end{equation*}
that is, $\mat{Y}_{p}^{\star}$ with an additional zero column appended. Since $\mat{Y}_{p}^{\star}$ is local minimizer for rank $p$, we have that $\mat{Y}_{0}$ is a critical point for the problem of rank $p+1$. As $\mat{Y}_{p}^{\star}$ is not the sought solution to \eqref{eq:formulation-on-psd}, this means that $\mat{Y}_{0}$ is a saddle point for the problem of rank $p+1$. Therefore, by virtue of the second order KKT optimality conditions, there must exists a descent direction $\mat{Z}\in\mathbb{R}^{n\times p}$ such that
\begin{equation*}
	\frac{1}{2}\trace(\mat{Z}^{T} D\grad\,f(\mat{Y}_{0})[\mat{Z}]) \leq 0.
\end{equation*}

To escape from the saddle point, we can thus exploit the following descent direction
\begin{equation*}
	\mat{Z} = [\mat{0}^{n\times p}|\vec{v}],
\end{equation*}
where $\vec{v}$ is the eigenvector associated to the smallest algebraic eigenvalue of 
\begin{equation}\label{eq:convex-gradient-at-local-min}
	\mat{S}_{Y} = \nabla_{\mat{X}}f(\mat{Y}_{p}^{\star}\mat{Y}_{p}^{\star T}),
\end{equation}
and where  $\nabla_{\mat{X}}f(\mat{Y}\mat{Y}^{T})$ is the Euclidean gradient of the convex cost function $f(\mat{X})$ evaluated at $\mat{Y}\mat{Y}^{T}$. As we have
\[
	\grad\,f(\mat{Y}_{0})=\nabla_{\mat{X}}f(\mat{Y}_{0}\mat{Y}^{T}_{0})\mat{Y}_{0},
	\] 
the proposed direction satisfies
\begin{equation*}
	\frac{1}{2}\trace(\mat{Z}^{T} D\grad\,f(\mat{Y}_{0})[\mat{Z}]) = \vec{v}^{T}\mat{S}_{Y}\vec{v} \leq 0.
\end{equation*}

The descent direction is exploited by performing a single line-search step using the Armijo rule. The resulting iterate is then used as the initial condition for the optimization algorithm that will solve the problem of rank $p+1$.

The procedure stops at the latest when $p=n$. However, in the setting of interest, problem \eqref{eq:formulation-on-psd} presents a low-rank solution with $r\ll n$. We thus, expect the algorithm to stop much before reaching $p=n$. 

For the proposed strategy it is important to reach a local minimum of the cost function as long as $p < r$. Although in theory convergence to saddle points cannot be excluded for gradient descent algorithms, the issue is not harmful in practice as saddle point are generally unstable from a numerical point of view. 

\section{DISCUSSION}
We propose both a gradient descent and a trust-region algorithm for solving the fixed-rank Euclidean distance matrix completion problem. The numerical cost per iteration for the gradient descent algorithm is $O(dp + np)$ versus $O(dp + np + np^2 + p^3)$ for the trust region algorithm. Although the gradient descent algorithm has a smaller computational cost per iteration, the number of iterations required to reach convergence is higher than for the trust-region algorithm.

We thus recommend the trust-region algorithm when a high optimization accuracy is required or when the observation noise is small. The gradient descent approach should be preferred for very large problems where the observation noise is high. In this setting, one is usually not interested in a solution of high-accuracy, since it generally compromises the generalization performance.

\section{NUMERICAL EXPERIMENTS}\label{sec:experiments}
In this section, we evaluate the performance of the proposed algorithms on benchmarks. A MATLAB implementation is available from the first author's webpage.\footnote{\url{http://www.montefiore.ulg.ac.be/~mishra}}

\subsection{A visual example}
This example is adapted from \cite{chu03a}. Consider $n=121$ data points arranged in a $3$-dimensional helix structure defined by
\begin{equation*}
	(x,y,z) = (4\cos(3t),4\sin(3t),2t),\quad 0\leq t \leq 2\pi.
\end{equation*}

\begin{figure}[ht!]
      \centering
      \includegraphics[width=0.85\textwidth]{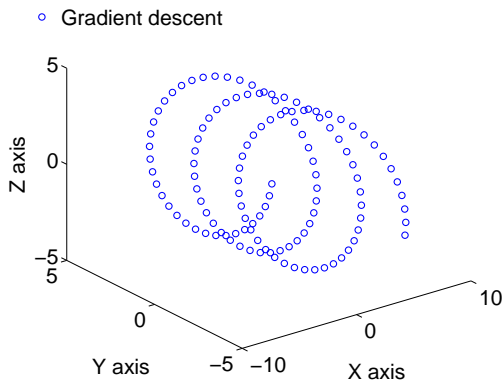}
      \caption{The proposed algorithms correctly recover the 3D helix structure form $15\%$ of the complete set of pairwise distances.}
      \label{fig:helix_gd}
   \end{figure}

After computing the distance matrix of between these points, we randomly remove $85\%$ of the distances uniformly and at random to generate a dissimilarity matrix $\widetilde{\mat{D}}$. From $15\%$ of distances the goal is to reconstruct the helix structure. We run the algorithms with the rank incremental strategy discussed in Section (\ref{sec:strategy-rank}). Both algorithms recover correctly the helix structure (Figure \ref{fig:helix_gd}). We only display the results for gradient descent as it coincides with the results of the trust-region algorithm. The algorithms stop when the relative or absolute variation of the cost function drops below $10^{-5}$.

\subsection{Trust-region versus gradient descent}
\begin{figure*}[htp]
\centering
\subfigure[ ]{
\includegraphics[scale=.26]{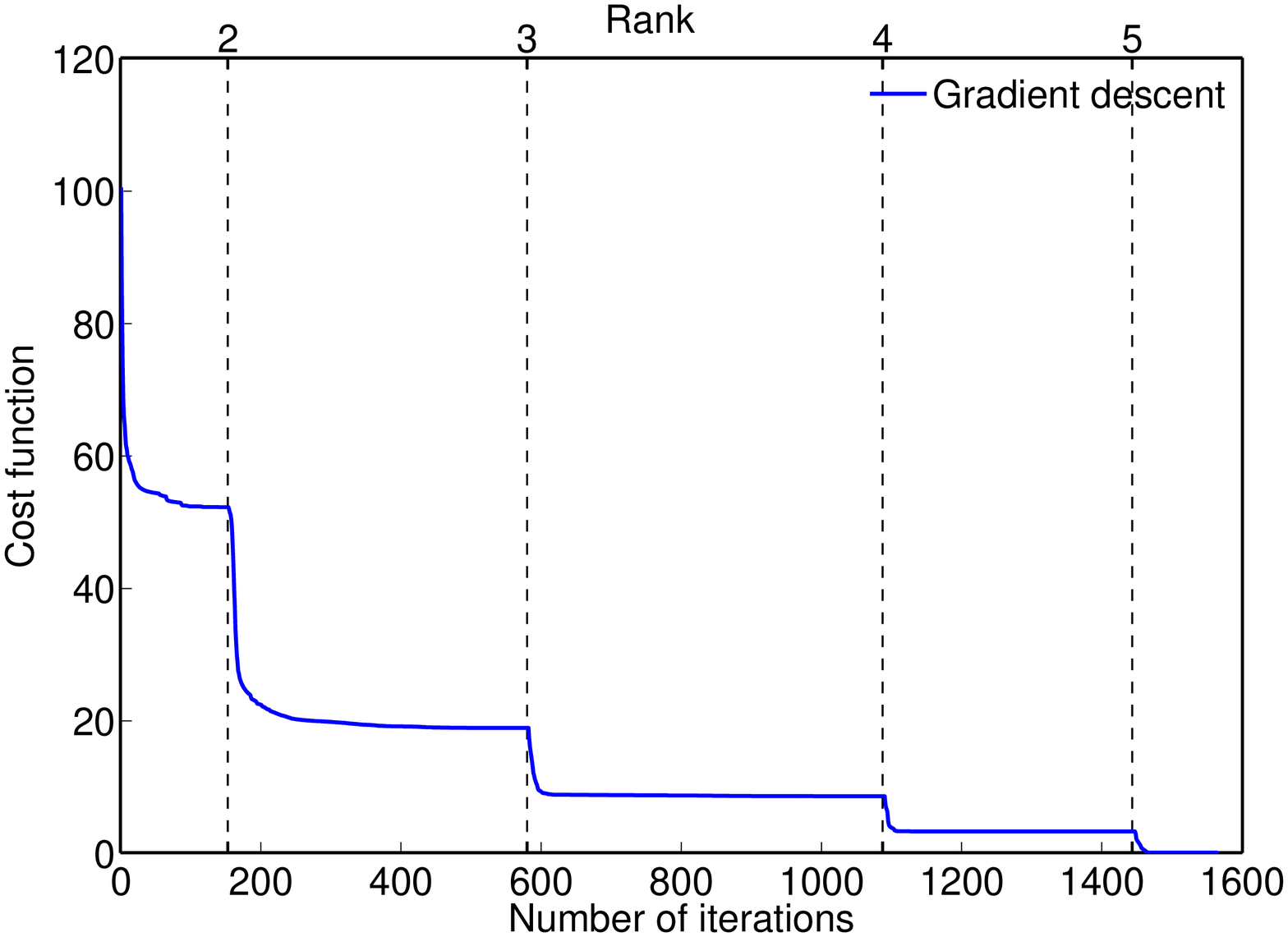}
\label{fig:convergence_gd}
}
\subfigure[ ]{
\includegraphics[scale=.26]{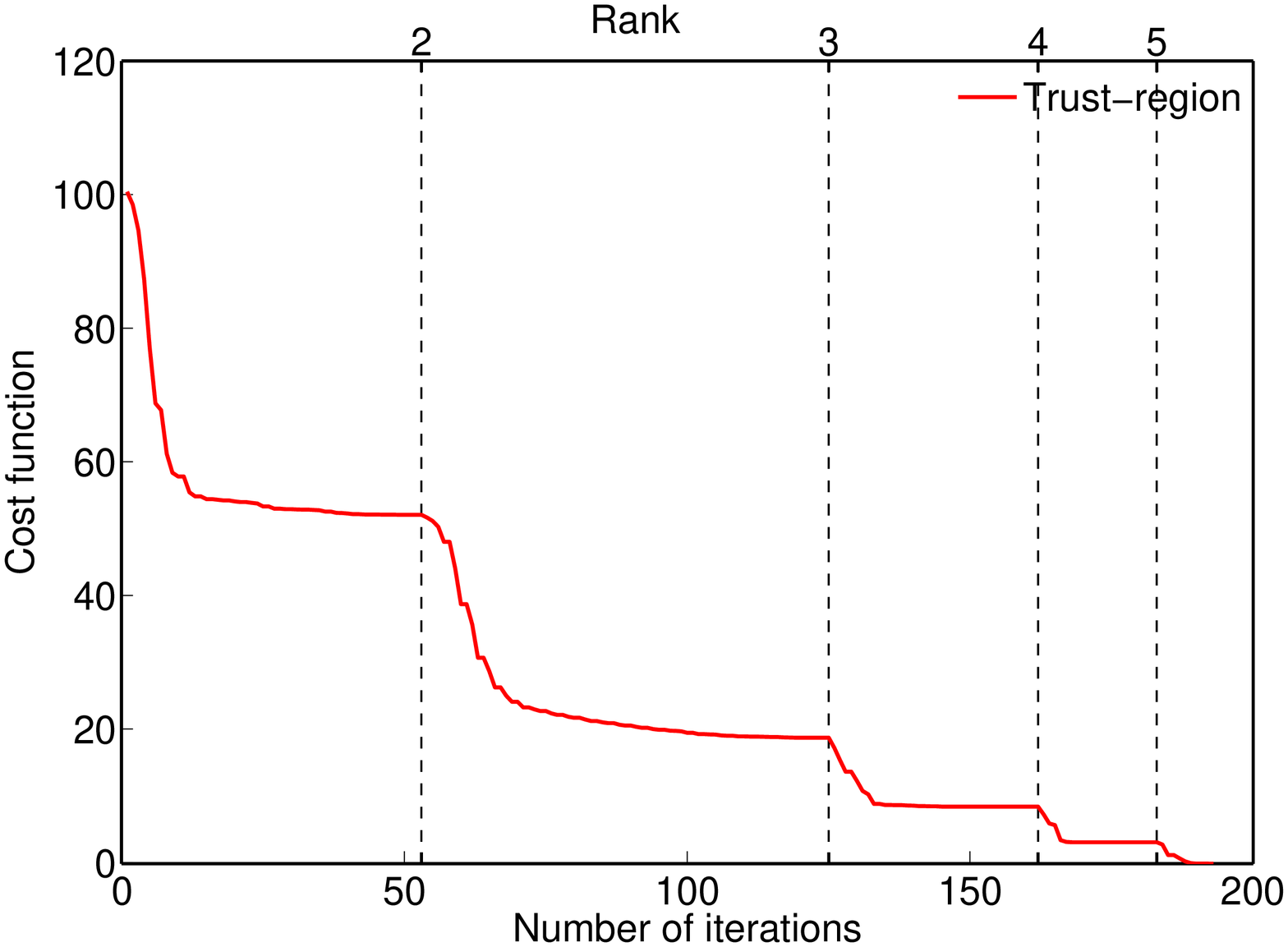}
\label{fig:convergence_tr}
}
\caption{Convergence plots of the algorithms.}
\end{figure*}

To compare the two versions of the algorithm, we generate a random distance matrix
\begin{equation}\label{eq:model-random}
	\mat{D}^{\star}= \kappa(\mat{Y}^{\star}\mat{Y}^{\star T}),
\end{equation}
where $\mat{Y}^{\star}\in\mathbb{R}^{500\times 3}$ has entries distributed according to gaussian distribution with zero mean and unit standard deviation. The fraction of unknown distances is fixed at $85\%$. We run the algorithms without knowing the embedding dimension. The algorithms are stopped when the relative or absolute variation of the cost function drops below $10^{-6}$. The objective function is plotted against the number of iterations (Figures \ref{fig:convergence_gd} and \ref{fig:convergence_tr}). Both algorithms recover the correct configuration and dimensionality. \change{The trust-region algorithm converges in $15.0$ seconds and $193$ iterations, whereas the gradient descent algorithm converges in $19.6$ seconds and $1565$ iterations.} Observe the monotonic convergence of both algorithms to the sought solution.

\subsection{Scaling test}
We now evaluate our algorithms on larger random data sets. We vary the problem size $n$ from $1000$ to $10000$. For each $n$, we generate a random distance matrix according to \eqref{eq:model-random} with $\mat{Y}^{\star}\in\mathbb{R}^{n\times 4}$. We sample $0.1$ fraction of the total amount of distances and the algorithms are run by fixing the embedding dimension, $p=4$. Results are averaged over $10$ runs. The test has been performed on a single core Intel L5420 2.5 GHz with 5GB of RAM.

\begin{figure*}[ht]
\centering
\subfigure[$\ $]{
\includegraphics[scale=.35]{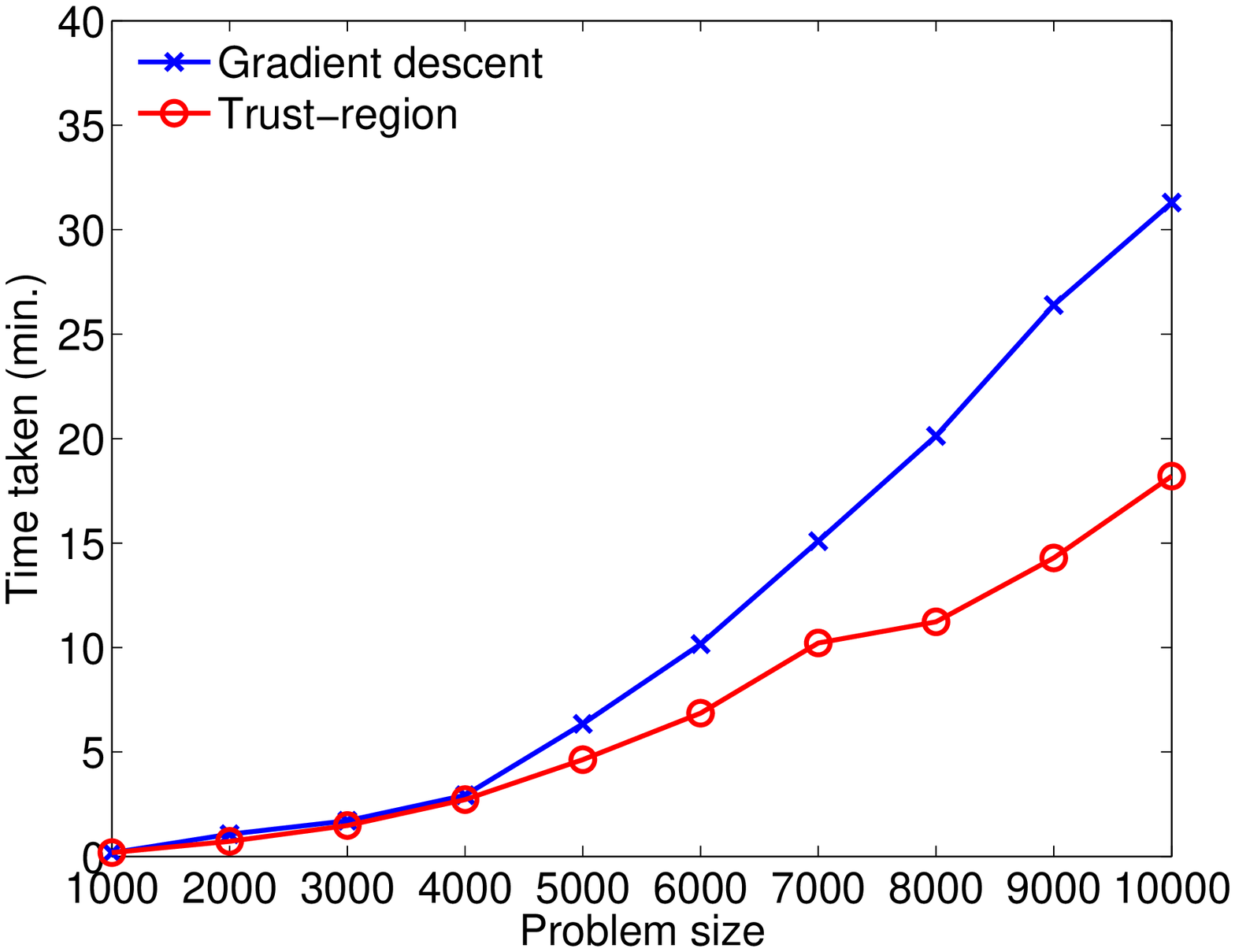} 
\label{fig:time_taken_test}    
}
\subfigure[$\ $]{
\includegraphics[scale=.35]{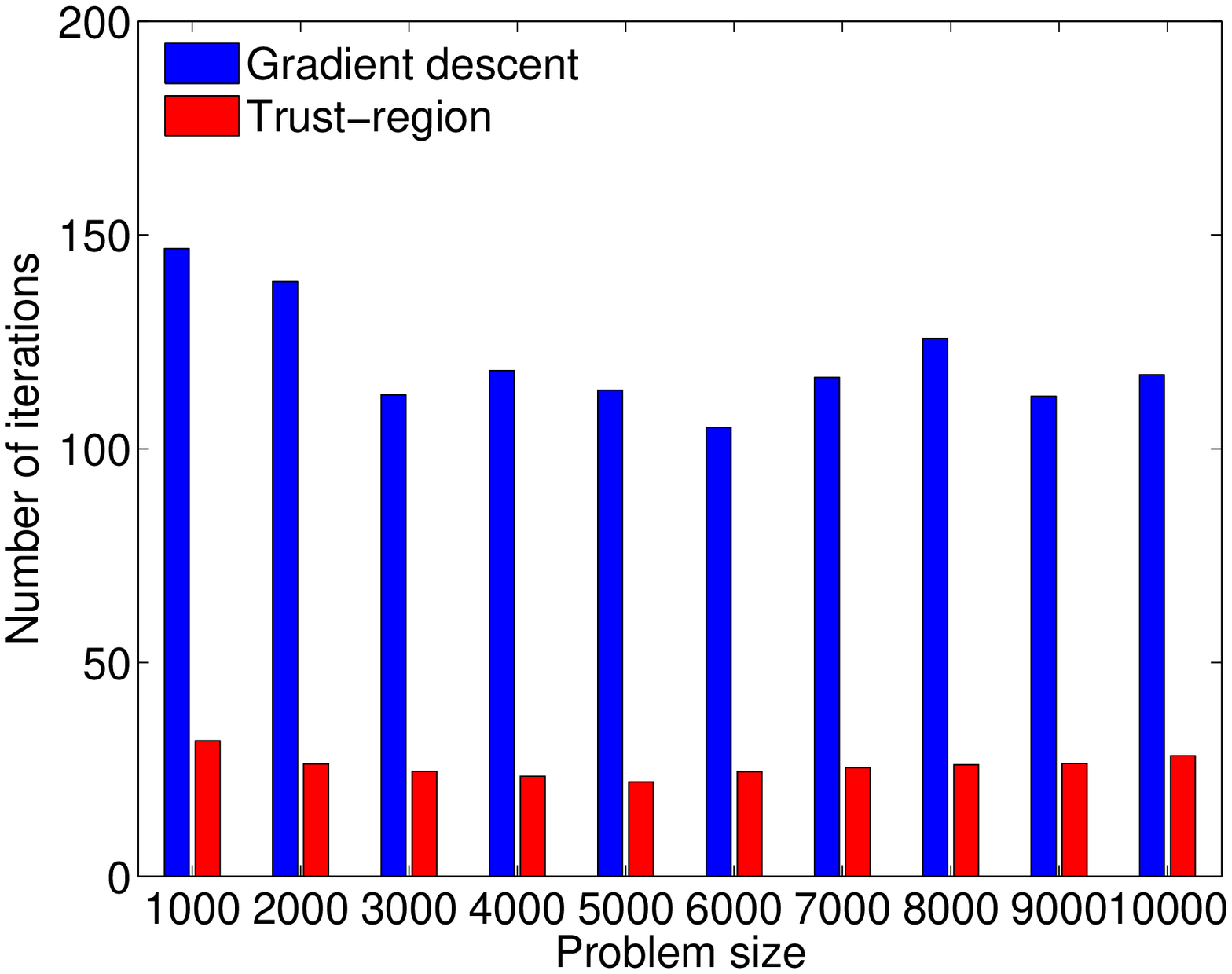}
\label{fig:iteration_test}
}
\caption{Analysis of the proposed algorithms on randomly generated datasets.}
\end{figure*}

The time taken and number of iterations required to reach convergence is reported at Figure \ref{fig:time_taken_test} and \ref{fig:iteration_test} respectively. For instance, for $n=10000$ the number of known distances is about $5$ millions ($10\%$ of $50$ million total entries). The gradient descent algorithm takes about $120$ iterations and 31 minutes, while the trust region algorithm solves the problem in $30$ iterations and $18$ minutes.

\section{CONCLUSION}\label{sec:conclusion}
In this paper two efficient numerical optimization algorithms have been presented for the distance matrix completion problem. In particular, the algorithms do not require any prior notion about the embedding and can potentially handle very large data sets. The proposed algorithms stem from a geometric view of the problem formulation. This interpretation as a manifold-based optimization problem considerably reduced the computational burden. At the same we were able to devise a superlinearly converging scheme namely, the trust-region algorithm in addition to the linearly convergent gradient descent algorithm. The numerical experiments that have been performed, are very encouraging on various parameters. 

\section{ACKNOWLEDGMENTS}
This paper presents research results of the Belgian Network DYSCO (Dynamical Systems, Control, and Optimization), funded by the Interuniversity Attraction Poles Programme, initiated by the Belgian State, Science Policy Office. The scientific responsibility rests with its authors. Gilles Meyer is supported as an FRS-FNRS research fellow (Belgian Fund for Scientific Research).

\bibliography{arXiv_MMS_cdc2011}
\bibliographystyle{IEEEtran}

\end{document}